\documentclass{amsart}

\usepackage{amsfonts}
\usepackage{amssymb}
\usepackage{amsmath}
\usepackage{graphicx}
\usepackage{color}

\def\({\left(}
\def\){\right)}

\numberwithin{equation}{section}

\newcommand{\cC}{\mathcal{C}}
\newcommand{\cD}{\mathcal{D}}
\newcommand{\cQ}{\mathcal{Q}}
\newcommand{\cS}{\mathcal{S}}
\newcommand{\cW}{\mathcal{W}}
\newcommand{\bN}{\mathbb{N}}
\newcommand{\bP}{\mathbb{P}}
\newcommand{\bQ}{\mathbb{Q}}
\newcommand{\bR}{\mathbb{R}}
\newcommand{\bZ}{\mathbb{Z}}
\newcommand{\fa}{\mathfrak{a}}
\newcommand{\fb}{\mathfrak{b}}
\newcommand{\SL}{\mathrm{SL}}
\newcommand{\ua}{\underline{a}}
\newcommand{\ub}{\underline{b}}
\newcommand{\uc}{\underline{c}}
\newcommand{\ud}{\underline{d}}
\newcommand{\uu}{\underline{u}}
\newcommand{\ur}{\underline{r}}
\newcommand{\uv}{\underline{v}}

\newcommand{\urj}{\underline{r_j}}
\newcommand{\uvj}{\underline{v_j}}

\newtheorem{theorem}{Theorem}[section]
\newtheorem{definition}[theorem]{Definition}
\newtheorem{proposition}[theorem]{Proposition}
\newtheorem{lemma}[theorem]{Lemma}
\newtheorem{remark}[theorem]{Remark}
\newtheorem{corollary}[theorem]{Corollary}

\makeatletter
\@namedef{subjclassname@2010}{%
  \textup{2010} Mathematics Subject Classification}
\makeatother

\begin{document}

\title{On Roots of Quadratic Congruences}

\author{Hieu T. Ngo}
\address{School of Applied Mathematics and Informatics, Hanoi University of Science and Technology, Hanoi, Vietnam}
\email{hieu.ngotrung@hust.edu.vn}

\subjclass[2010]{11M20} 
\keywords{Binary quadratic forms, congruence subgroups, Kloosterman sums, Poincar\'e series.}

\begin{abstract} 
The equidistribution of roots of quadratic congruences with prime moduli depends crucially upon effective bounds for a special Weyl linear form. Duke, Friedlander and Iwaniec discovered a strong estimate for this Weyl linear form when the quadratic polynomial has negative discriminant. T\'oth established an analogous but weaker bound when the quadratic polynomial has positive discriminant. We obtain a stronger estimate for the Weyl linear form for quadratics of positive discriminants.
\end{abstract}

\maketitle

\section{Introduction}


Let $f(X)=\alpha X^2 + \beta X + \gamma$ be an irreducible quadratic polynomial with integer coefficients of discriminant $\Delta=\beta^2-4\alpha\gamma$. The distribution of the roots of the congruence equation
\begin{equation*}
f(x) \equiv 0 \, (\rm{mod}\, n)
\end{equation*}
as the integer modulus $n$ varies is a classical topic in number theory.
The `quadratic harmonic'
$$
\rho_h(n)=\sum_{\substack{\nu \,{\rm mod} \, n \\ f(\nu)\equiv 0 \, (\rm{mod}\, n)}} e\(\frac{h\nu}{n}\),
$$
where $h$ is an integer, plays a fundamental role in the analysis of roots of quadratic congruences.
One would expect that there is nontrivial cancellation for the discrete Weyl linear form
\begin{equation*}
\cW_h(x,N) = \sum_{\substack{x < n < 2x \\ n \equiv 0 \, (\rm{mod}\, N)}} \rho_h(n)
\end{equation*}
for every nonzero integer $h$. For an irreducible quadratic polynomial $f(X)$ with integer coefficients and for fixed $h \in \bZ \backslash \{0\}$, Hooley \cite[Theorem 1]{Hooley1} proved
\begin{equation}\label{eq:Hooley-bound}
\cW_h(x,1) \ll_h x^{\frac34} \(\log x\)^2 ,
\end{equation}
establishing equidistribution of roots of quadratic congruences when the moduli are integers. Hooley generalized in \cite{Hooley2} that roots of polynomial congruences modulo integers are equidistributed when the defining polynomial is irreducible of degree at least two.
Bykovskii \cite[Theorem 4]{By} and Hejhal \cite[page 332]{Hejhal} employed spectral methods to improve the exponent $\frac34$ in \eqref{eq:Hooley-bound} to $\frac23 + \epsilon$. 

It is natural to ask whether roots of quadratic congruences modulo prime numbers are equidistributed, and it turns out that this question is much more challenging than the original problem for the integers. Using sieve methods and spectral analysis of Poincar\'e series, Duke, Friedlander and Iwaniec \cite{DFI1} established the startling theorem asserting that roots of quadratic congruences with prime moduli are equidistributed when the irreducible quadratic polynomial defining the congruences has negative discriminant. T\'oth \cite{Toth} introduced several technical innovations and established the analogous theorem for an irreducible quadratic polynomial of positive discriminant. In both of the papers \cite{DFI1} and \cite{Toth}, the crux of the proof of the equidistribution theorem is a sufficiently strong estimate for the discrete Weyl linear form. More precisely, Duke, Friedlander and Iwaniec \cite[Proposition 1]{DFI1} showed that, for an irreducible quadratic polynomial of negative discriminant,  
\begin{equation}\label{eq:Weyl-form-DFI}
\cW_h(x,N) \ll_\epsilon \gcd(h,N)^{\frac{1}{20}} \(\frac{N^2}{x}\)^{\frac{1}{20}} \(\frac{x}{N}\)^{1+\epsilon}. 
\end{equation}
For an irreducible quadratic polynomial of positive discriminant, T\'oth \cite[p.~738]{Toth} showed that for all sufficiently large natural numbers $L$ one has
\begin{equation}\label{eq:Weyl-form-Toth}
\cW_h(x,N) \ll_{h} \(\frac{N^2}{x}\)^{\frac{1}{4L}} \(\frac{x}{N}\)^{1+\frac{1}{L^2}}. 
\end{equation} 
The bounds \eqref{eq:Weyl-form-DFI} and \eqref{eq:Weyl-form-Toth} provide nontrivial cancellation in the important range $N^2 = o(x)$.
One observes that there is a gap of quality between the two estimates: the factor $\frac{1}{L^2}$ makes the bound \eqref{eq:Weyl-form-Toth} considerably weaker than the bound \eqref{eq:Weyl-form-DFI}.

The main purpose of this paper is to establish a stronger bound for the discrete Weyl linear form in the positive-discriminant case. 

\begin{theorem}\label{thm:discrete-linear-form-estimate}
Let $f(X)$ be an irreducible quadratic polynomial with integer coefficients whose discriminant is positive; let $h$ be a nonzero integer. One has
\begin{equation}\label{eq:discrete-linear-form-estimate}
\cW_{h}(x,N) \ll_\epsilon  \(x^{\frac{12}{13}} N^{-\frac{11}{13}} h^{\frac{1}{13}}  +  h\) x^\epsilon . 
\end{equation}
\end{theorem}

\begin{remark}
{\rm An equivalent formulation of \eqref{eq:discrete-linear-form-estimate} is to assert that
\begin{equation}\label{eq:Weyl-form-H}
\cW_{h}(x,N) \ll_\epsilon  
	\(\frac{N^2}{x}\)^{\frac{1}{13}} \(\frac{x}{N}\)^{1+\epsilon} h^{\frac{1}{13}}  +  h x^\epsilon . 
\end{equation}
In this form, the estimate \eqref{eq:Weyl-form-H} provides nontrivial cancellation in the range $N^2=o(x)$ and improves upon \eqref{eq:Weyl-form-Toth} for quadratics of positive discriminants. The numerical strength of \eqref{eq:Weyl-form-H} for the positive-discriminant case is (slightly) better than that of \eqref{eq:Weyl-form-DFI} for the negative-discriminant case when $h$ is fixed because the saving factor $\frac{N^2}{x}$ has exponent $\frac{1}{13}>\frac{1}{20}$. Our method seems to work for the negative-discriminant case. We leave the consideration of this case together with further consequences to our next work.}
\end{remark}

Although our proof proceeds along a similar line to that of T\'oth, there are several differences that are both subtle and necessary to strengthen the estimate of the Weyl linear form. We put more care into the smoothing step, letting the weight function's parameters participate in the spectral analysis. We use a slightly different normalization in the transition from roots of quadratic congruences to Poincar\'e series which paves the way for Kloosterman sums on Hecke congruence subgroups. The key new ingredient is a strong bound on sums of Kloosterman sums which was essentially due to Pitt and used in a very different context. The proof of Pitt's estimate involves a nontrivial analysis of the exceptional spectrum.

The remainder of this paper is organized as follows. In Section \ref{section:preliminary}, we collect facts and results concerning Kloosterman sums for congruence groups. 
The analysis of Weyl linear forms is carried out in Section \ref{section:Weyl-linear-forms}, where the proof of Theorem \ref{thm:discrete-linear-form-estimate} is presented. We first smooth the discrete Weyl linear form in Section \ref{subsection:Weyl-linear-form}, next express the smooth Weyl linear forms as Poincar\'e series in Sections \ref{subsection:binary-quadratic-form} -- \ref{subsection:partition-of-unity}, then spectrally analyse the Poincar\'e series in Section \ref{subsection:Poincare-Kloosterman}, and finally complete the proofs of some technical lemmas in Section \ref{subsection:auxiliary-estimate}. 
 
\section{Preliminaries}\label{section:preliminary}




\subsection{Kloosterman sums for Hecke congruence subgroups}

We review congruence groups and Kloosterman sums, mostly following the presentation and notation of Iwaniec in his books \cite[Chapters 2 and 4]{Iw1} and \cite[Chapter 2]{Iw2}.

Let $\Lambda\subset\SL_2(\bZ)$ be a subgroup of finite index; suppose that $-\begin{pmatrix} 1 & 0 \\ 0 & 1 \end{pmatrix} \in \Lambda$. An element of $\Lambda$ is called \emph{parabolic} if it has a unique fixed point in $\bP^1(\bR) = \bR\cup \{\infty\}$. A point $\fa \in \bP^1(\bQ) = \bQ\cup \{\infty\}$ is called a \emph{cusp} of $\Lambda$ if there exists a parabolic element of $\Lambda$ which fixes $\fa$; denote by $\Lambda_\fa$ the stabilizer of $\fa$. 
Two cusps $\fa$ and $\fb$ of $\Lambda$ are \emph{equivalent} if the two orbits $\Lambda \fa$ and $\Lambda \fb$ are the same.
For each cusp $\fa$ of $\Lambda$, a \emph{scaling matrix for $\fa$} is an element $\sigma_\fa \in \SL_2(\bR)$ which satisfies $\sigma_\fa(\infty)=\fa$ and $\sigma_\fa^{-1}\Lambda_\fa \sigma_\fa = \pm \begin{pmatrix} 1 & \bZ \\ 0 & 1 \end{pmatrix} =: P_0$; any two scaling matrices $\sigma_\fa$ and $\sigma'_\fa$ for the same cusp $\fa$ satisfy the relation $\sigma'_\fa=\sigma_\fa\cdot \begin{pmatrix} t & r \\ 0 & t^{-1} \end{pmatrix}$ for some $r\in\bR$ and $t \in \{\pm 1\}$. We associate to any two cusps $\fa$ and $\fb$ the set
\begin{equation}\label{eq:cusps-moduli}
\cC_{\fa \fb}(\Lambda) = \{ c \in \bR_{>0}: \textrm{there exist } a,b,d \in \bR \textrm{ so that } \begin{pmatrix} a & b \\ c & d \end{pmatrix} \in \sigma_\fa^{-1} \Lambda \sigma_\fb \};
\end{equation}
the set $\cC_{\fa \fb}(\Lambda)$ depends only on the cusps $\fa$ and $\fb$ and is independent of the choices of the scaling matrices $\sigma_\fa$ and $\sigma_\fb$. For every $c\in\cC_{\fa \fb}(\Lambda)$, define
\begin{align*}
\cD_{\sigma_\fa \sigma_\fb}(c) &= 
	\{ d\in\bR: 0<d\leq c, \textrm{there exist } a,b \in \bR \textrm{ so that } 
	\begin{pmatrix} a & b \\ c & d \end{pmatrix} \in \sigma_\fa^{-1} \Lambda \sigma_\fb \} ; 
\end{align*}
the set $\cD_{\sigma_\fa \sigma_\fb}(c)$ depends on the choices of the scaling matrices $\sigma_\fa$ and $\sigma_\fb$. Observe that for each $d \in \cD_{\sigma_\fa \sigma_\fb}(c)$, the value $a \, ({\rm mod \, } c)$ such that there exists $b\in\bR$ with
$\begin{pmatrix} a & b \\ c & d \end{pmatrix} \in \sigma_\fa^{-1} \Lambda \sigma_\fb$ is uniquely determined. We are in a position to define Kloosterman sums.

\begin{definition}\label{defn:Kloosterman-sum}
Let $\fa$ and $\fb$ be two cusps of $\Lambda$; let $\sigma_\fa$ and $\sigma_\fb$ be two scaling matrices for $\fa$ and $\fb$ respectively. Let $m,n\in\bZ$ and $c\in \cC_{\fa \fb}(\Lambda)$. The \emph{Kloosterman sum associated to $\sigma_\fa$ and $\sigma_\fb$ with modulus $c$ and frequencies $m,n$} is
\begin{align}
\cS_{\sigma_\fa \sigma_\fb}(m,n;c) 
	&= \sum_{\begin{pmatrix} a & b \\ c & d \end{pmatrix} \in 
	P_0 \backslash \sigma_\fa^{-1} \Lambda \sigma_\fb / P_0}
	 e\(\frac{ma+nd}{c}\)
	= \sum_{d \in \cD_{\sigma_\fa \sigma_\fb}(c)}
	e\(\frac{ma+nd}{c}\)    \label{eq:Kloosterman} \\
	&= \sum_{\gamma \in 
	\Lambda_\fa \backslash \Lambda / \Lambda_\fb, 
	\begin{pmatrix} a & b \\ c & d \end{pmatrix} := \sigma_\fa^{-1} \gamma \sigma_\fb}
	e\(\frac{ma+nd}{c}\)  .  \label{eq:Kloosterman1} 
\end{align}
\end{definition}

\begin{definition}
Let $\eta\in\SL_2(\bZ)$ and put $\fa= \eta\cdot\infty \in \bP^1(\bQ)$. Let $\Lambda_\fa$ be the stabilizer in $\Lambda$ of $\fa$.
Observe that the group $\eta^{-1} \Lambda_\fa \eta$ depends only on the coset $\Lambda\eta$ in $\Lambda\backslash\SL_2(\bZ)$; for if $\eta'=\gamma\eta$ with $\gamma\in\Lambda$ and $\fa'= \eta'\cdot\infty$, then $\Lambda_{\fa'} = \gamma \Lambda_\fa \gamma^{-1}$ and hence
$\eta'^{-1} \Lambda_{\fa'} \eta' = \eta^{-1} \Lambda_\fa \eta$.
We call the smallest positive integer $w$ such that
$$
\eta^{-1} \Lambda_\fa \eta = \pm  \begin{pmatrix} 1 & w\bZ \\ 0 & 1 \end{pmatrix}
$$
the \emph{width at infinity of the coset  $\Lambda\eta$ in $\Lambda\backslash\SL_2(\bZ)$}, denoted by $w={\rm width}_\Lambda (\Lambda\eta)$. We have a well-defined map ${\rm width}_\Lambda : \Lambda\backslash\SL_2(\bZ) \to \bZ_{\geq 0}$.
\end{definition}

Let us now specialize to the Hecke congruence subgroup $\Lambda=\Gamma_0(q)$ of level $q$. Our goal is to describe the Kloosterman sums $\cS_{\sigma_\infty \sigma_\fa}(m,n;c)$ which are pertinent to roots of quadratic congruences. 
Explicit formulas for these Kloosterman sums are not easy to extract from the literature. Motohashi's papers \cite{Motohashi1} (for $\Gamma_0(q)$ with $q$ squarefree), \cite[Sections 12--15]{Motohashi2} and \cite{Motohashi3} (for $\Gamma_0(q)$ with a general modulus $q$) contain computations relevant to our discussion. See also \cite{Blomer, Watt}.

The following description of cusps of Hecke congruence subgroups is standard (cf.~\cite[Section 2.5]{Iw1}, \cite[Section 2.3]{Iw2}, \cite[Proposition 1.43]{Shimura}).

\begin{lemma}\label{lem:cusps-representatives-width}
A complete set of cusps which are inequivalent modulo $\Gamma_0(q)$, i.e.~a complete set of representatives for $\Gamma_0(q)\backslash \bP^1(\bQ)$ is
\begin{equation}\label{eq:cusps}
\Big\{ \frac{\mu}{\nu}: \mu,\nu \in \bZ_{>0}, \nu | q, \gcd(\mu,\nu)=1, \mu \textrm{ mod } \, \gcd(\nu, q/\nu) \Big\}.
\end{equation}
The number of inequivalent cusps for $\Gamma_0(q)$ is $\sum_{\nu > 0, \nu | q}\varphi(\gcd(\nu, q/\nu))$.
\end{lemma}

The definition of Kloosterman sums depends on the choice of a scaling matrix $\sigma_\fa$ for each cusp $\fa$ of $\Gamma_0(q)$. We choose $\sigma_\infty = \begin{pmatrix} 1 & 0 \\ 0 & 1 \end{pmatrix}$.
By Lemma \ref{lem:cusps-representatives-width}, suppose that $\fa=\frac{\mu}{\nu}\in\bQ$ is a cusp for $\Gamma_0(q)$ with $\mu,\nu\in\bN$ satisfying $\gcd(\mu,\nu)=1$ and $\nu | q$. Put $q=\nu q'$, $\iota=\gcd(\nu,q')$, $\nu=\iota\nu'$, and $q'=\iota q''$; we have $\gcd(\nu',q'')=1$. 
Choose an integer $\overline{\mu} \in \bZ$ such that $\mu\overline{\mu} \equiv 1  \, (\rm{mod}\, \nu)$.
Set 
$$
\eta_\fa =  \begin{pmatrix} \mu & \frac{\mu\overline{\mu}-1}{\nu} \\ \nu & \overline{\mu} \end{pmatrix} ,
\quad \tau_\fa =  \begin{pmatrix} \sqrt{q''} & 0 \\ 0 & 1 / \sqrt{q''} \end{pmatrix} ,
$$
so that $\eta_\fa\cdot\infty=\fa$.
Then $\sigma_\fa = \eta_\fa \tau_\fa$ is a scaling matrix for $\fa$.
The stabilizer in $\Gamma_0(q)$ of $\fa$ is
$$
\Gamma_0(q)_\fa	
	= \bigg\{
	\begin{pmatrix}
		1+l\fa & -l\fa^2 \\
		l & 1-l\fa 
	\end{pmatrix} :
	l \equiv 0 \,\,{\rm mod}\, [q,\nu^2]
	\bigg\}
$$
where $[q,\nu^2]$ denotes the least common multiple of $q$ and $\nu^2$.
We are in a position to describe the set of Kloosterman moduli $\cC_{\infty \fa}(\Gamma_0(q))$.

\begin{lemma}\label{lem:Kloosterman-moduli}
Suppose that $\eta \in \SL_2(\bZ)$ and that $\fa=\eta\cdot\infty=\frac{\mu}{\nu}\in\bQ$ is a cusp for $\Gamma_0(q)$ with $\gcd(\mu,\nu)=1$ and $\nu | q$. Put $q=\nu q'$, $\iota=\gcd(\nu,q')$, $\nu=\iota\nu'$, and $q'=\iota q''$. 
\begin{enumerate}
\item One has 
$${\rm width}_{\Gamma_0(q)}(\Gamma_0(q)\eta)=\frac{q}{\gcd(q,\nu^2)}=q''.$$
\item Let
$$
\cQ(\nu,q') = \{ \nu c: c\in\bN, \gcd(c,q') = 1 \}.
$$
Then $\cC_{\infty \fa}(\Gamma_0(q)) = \sqrt{q''} \cQ(\nu,q')$.
\end{enumerate}
\end{lemma}

\begin{proof}
The width in $(1)$ was calculated in \cite[Section 13]{Motohashi2}. The statement $(2)$ can be deduced from an easy but tedious calculation.
\end{proof}

\subsection{Sums of Kloosterman sums}
We shall need a good bound on sums of Kloosterman sums.
The following estimate is essentially due to Pitt \cite{Pitt}. 

\begin{theorem}\label{thm:Pitt-estimate}
Let $V(c,\kappa)$ be a smooth and compactly supported function on $C < c < 2C$ and $K < \kappa < 2K$. Suppose that there is a real number $Y\geq 1$ such that for all $I,J\in\bN$, the derivatives of $V(c,\kappa)$ is bounded by
$$
\( \frac{d}{dc} \)^I \( \frac{d}{d\kappa} \)^J  V(c,\kappa) 
	\ll  \frac{Y^{I+J}}{C^I K^J}  (YCK)^\epsilon .
$$
Then
\begin{align}
&\sum_{\kappa\in\bZ} \,\, \sum_{c\in\cC_{\infty\fa}(\Gamma_0(q))/\sqrt{q''}} 
	\frac{S_{\sigma_\infty \sigma_\fa}(h,\kappa;c\sqrt{q''})}{c} V(c,\kappa) \nonumber \\
&\ll K^{\frac12} \bigg\{ 
	 \( K^\frac14 q^{-\frac14} + 1 \) C^\frac12 Y^\frac32 \gcd(h,q)^\frac14 
	+
	q^{\frac12}\( Y^\frac34 + K^\frac12 q^{-\frac12} \) Y^\frac74    
	 \bigg\} 
	(qYKC)^\epsilon  .   \label{eq:Pitt-theorem} 
\end{align}
\end{theorem}

Theorem \ref{thm:Pitt-estimate} is a prototype of the interplay between Kloosterman sums and automorphic forms for Hecke congruence subgroups, the foundation of which was built by Deshouillers and Iwaniec \cite{DI}.
In fact, Pitt proved his bound \cite[Theorem 1.6]{Pitt} for a very similar Kloosterman sum of the congruence subgroup $\Gamma_0(q)$.
The proof of Pitt's theorem employed Kuznetsov trace formula \cite{Kuznetsov} and spectral large sieve inequalities for $\Gamma_0(q)$, including considerations of exceptional eigenvalues. 
Exactly the same arguments can be used to establish Theorem \ref{thm:Pitt-estimate}. 






\section{Analysis of Weyl linear forms}\label{section:Weyl-linear-forms}


\subsection{Weyl linear forms}\label{subsection:Weyl-linear-form}

We are interested in exhibiting cancellation for the discrete Weyl linear form
\begin{equation*}
\cW_h(x,N) = \sum_{\substack{x\leq n\leq 2x \\ n \equiv 0 \, (\rm{mod}\, N)}} \rho_h(n).
\end{equation*}
A standard approach is to approximate $\cW_h(x;N)$ by a smooth Weyl linear form. Let $1 < Y_1 < x$ be a parameter to be determined. Let $g:\bR\to [0,1]$ be a weight function such that $g$ is smooth and compactly supported in the interval $[x-\frac{x}{Y_1},2x+\frac{x}{Y_1}]$, that $g$ equals $1$ on the interval $[x,2x]$, and that its derivatives are bounded by $\| g^{(j)} \|_\infty \ll_j \(\frac{Y_1}{x}\)^j$ for all $j\in\bN$.
Define the smooth Weyl linear form
\begin{equation}\label{eq:smooth-linear-form}
\cW_{h,g}(x,N) = \sum_{n \equiv 0 \, (\rm{mod}\, N)} \rho_h(n) g(n).
\end{equation}

\begin{lemma}\label{lem:approximate-linear-form} We have
\begin{equation}\label{eq:approximate-linear-form}
|\cW_h(x,N) - \cW_{h,g}(x,N)| \ll \frac{x \tau(N) \log x}{NY_1} + x^\frac12
\end{equation}
where $\tau(N)$ denotes the number of divisors of $N$.
\end{lemma}

\begin{proof}
Using the fact that $|\rho_h(n)| \ll \tau(n)$, we have
\begin{align*}
|\cW_h(x,N) - \cW_{h,g}(x,N)| 
	&\ll  \sum_{\substack{x-\frac{x}{Y_1} \leq n \leq x \\ n \equiv 0 \, ({\rm mod} \, N)}} \tau(n)
		+ \sum_{\substack{2x \leq n \leq 2x + \frac{x}{Y_1}\\ n \equiv 0 \, ({\rm mod} \, N)}} \tau(n)	
\end{align*}
We now use the inequality $\tau(uv) \leq \tau(u)\tau(v)$ for $u,v\in\bN$ and the asymptotic
$$
\sum_{n\leq y}\tau(n) = y\log y + (2\gamma-1)y + O\(y^\frac12\),
$$
where $\gamma=\lim_{n\to+\infty}\(\log n + \sum_{k=1}^n\frac{1}{k}\)=0.5772\dots$ denotes the Euler's constant, to conclude the lemma.
\end{proof}

\subsection{Binary quadratic forms}\label{subsection:binary-quadratic-form}

We recall a classical correspondence between roots of quadratic congruences and binary quadratic forms, which has the effect of relating the smooth Weyl linear form $\cW_{h,g}(x,N)$ to automorphic functions. 
This connection was used in many works which study quadratic congruences (see \cite{Hooley1, By, Hejhal, DFI1, DFI2, Toth}).
We mainly follow T\'oth \cite{Toth}, who advocated the use of group-theoretical language in this context. However, our normalization is different from that of \cite{Toth}; we normalize in such a way that the congruence groups that arise are more appropriate for our study. The difference is subtle, but it paves the way for our spectral analysis.

For $\xi = \begin{pmatrix} a & b \\ c & d \end{pmatrix} \in \SL_2(\bR)$, we view the matrix coefficients $a,b,c,d$ as functions on $\SL_2(\bR)$ and denote these functions as $\ua(\xi), \ub(\xi), \uc(\xi), \ud(\xi)$ respectively.
We abbreviate a binary quadratic form as $uX^2 + rXY + vY^2=[u,r,v]$.
The group $\SL_2(\bR)$ acts on the set of real quadratics
$$\cQ_\bR = \{ [u,r,v]: u,r,v\in\bR \}$$
by linear changes of variables
$$ (\xi\cdot q)(X,Y) = q((X,Y)\xi) = q(aX+cY,bX+dY) \qquad \text{for } \xi  = \begin{pmatrix} a & b \\ c & d \end{pmatrix} \in \SL_2(\bR).$$
In other words, if $q=[u,r,v]$ with $u,r,v\in\bR$, by definition we have
\begin{equation}\label{eq:group-action}
\xi \cdot q = \xi \cdot [u,r,v] = [\uu(\xi),\ur(\xi),\uv(\xi)]
\end{equation}
where
\begin{align}
\uu(\xi) &=  q(a,b) = ua^2+rab+vb^2 ,  \label{eq:group-action-u} \\
\ur(\xi) &=   (2uac+2vbd) + r(ad+bc) , \label{eq:group-action-r} \\
\uv(\xi) &= q(c,d) = uc^2+rcd+vd^2. \label{eq:group-action-v}
\end{align}
Note that we have viewed the form coefficients as functions on the unimodular group $\SL_2(\bR)$.

Consider the quadratic polynomial $f(X)=\alpha X^2 + \beta X + \gamma$ with discriminant $\Delta=\beta^2-4\alpha\gamma$. 
The restriction of the action \eqref{eq:group-action} to the Hecke congruence subgroup 
\begin{equation}\label{eq:group}
\Gamma := \Gamma_0(\alpha) 
	= \bigg\{\begin{pmatrix} a & b \\ c & d \end{pmatrix} 
		\in \SL_2(\bZ) : c \equiv 0  \, (\rm{mod}\, \alpha) \bigg\}
\end{equation}
preserves the set of quadratics pertaining to $f$
\begin{equation}\label{eq:set-quadratics}
\cQ_f = \{ [u,r,v]: u,r,v\in\bZ, v \equiv 0 \, (\rm{mod}\, \alpha) , r \equiv \beta \, (\rm{mod}\, 2\alpha) , r^2-4uv=\Delta \}.
\end{equation}
In other words, if $\xi\in\Gamma$ and $q\in\cQ_f$, then $\xi\cdot q\in \cQ_f$.
As a consequence, given a quadratic $[u,r,v] \in \cQ_f$, we may view the form coefficients as functions $\uu(\xi),\ur(\xi),\uv(\xi)$ of $\xi\in\Gamma$.
The set of orbits $\Gamma\backslash\cQ_f$ is finite. Let us choose once and for all a complete set of representatives $\{q_j=[u_j,r_j,v_j]:1\leq j\leq h_\Delta\} \subset \cQ_f$ for the orbit set $\Gamma\backslash\cQ_f$.

If $1\leq j\leq h_\Delta$, we write $\Gamma^{(j)}=\{ \gamma \in \Gamma : \gamma \cdot q_j = q_j \}$ for the group of automorphs of $q_j$. The isotropy groups $\Gamma^{(j)}$ can be described in terms of solutions of Pell equations. Let $(\tau_0,\upsilon_0)$ be the fundamental solution of the Pell equation $\tau^2 - \Delta \upsilon^2=4$; in other words $(\tau_0,\upsilon_0)\in\bN^2$ satisfies 
$\tau_0^2 - \Delta \upsilon_0^2=4$ and if $(\tau,\upsilon)\in\bN^2$, $(\tau,\upsilon) \neq (\tau_0,\upsilon_0)$, $\tau^2 - \Delta \upsilon^2=4$, then 
$\tau_0+\sqrt{\Delta}\upsilon_0 < \tau+\sqrt{\Delta}\upsilon$.
We have (see 
\cite[Theorems 3.9 and 3.10]{Buell}) 
\begin{align}
\Gamma^{(j)} &= 
	\bigg\{  \begin{pmatrix} \frac{\tau+\upsilon r_j}{2} & -\upsilon u_j \\ 
		\upsilon v_j & \frac{\tau-\upsilon r_j}{2} \end{pmatrix} : 
		\tau,\upsilon \in \bZ, \tau^2 - \Delta \upsilon^2 = 4 \bigg\} \label{eq:isotropy-group1} \\
	&= \bigg\{ \pm \begin{pmatrix} \frac{\tau_0+\upsilon_0 r_j}{2} & -\upsilon_0 u_j \\ 
		\upsilon_0 v_j & \frac{\tau_0-\upsilon_0 r_j}{2} \end{pmatrix}^n : 
		n \in \bZ \bigg\} .  \label{eq:isotropy-group2}
\end{align}

We recall an identity which is due to Hooley \cite[Equation (27) page 291]{Hooley3}.
\begin{lemma}\label{lem:Hooley-identity}
For $\xi = \begin{pmatrix} a & b \\ c & d \end{pmatrix} \in \SL_2(\bR)$ and $[u,r,v]\in\cQ_f$, we have
\begin{equation}\label{eq:Hooley-identity}
\frac{\ur(\xi)}{\uv(\xi)} = \frac{2a}{c} - \frac{rc+2vd}{c \uv(\xi)}.
\end{equation}
\end{lemma}
\begin{proof}
The identity \eqref{eq:Hooley-identity} follows immediately from \eqref{eq:group-action-r}, \eqref{eq:group-action-v}, and the assumption $\det(\xi)=1$.
\end{proof}

Let us consider a congruence equation $f(\nu) \equiv 0 \, (\rm{mod}\, n)$. On completing the square $\eta=2\alpha\nu+\beta$, we obtain a bijection between the set of congruence roots
\begin{equation}\label{eq:set-congruence-roots}
Z_f(n) = \{ \nu \, (\rm{mod}\, n): f(\nu) \equiv 0 \, (\rm{mod}\, n)\}
\end{equation}
and the set
\begin{equation}\label{eq:set-congruence-roots-square}
Z'_f(n) = \{ \eta \, (\rm{mod}\, 2\alpha n):  
\eta \equiv \beta \, (\rm{mod}\, 2\alpha) , \eta^2 \equiv \Delta \, (\rm{mod}\, 4\alpha n)\}.
\end{equation}
By classical theory one connects the set of congruence roots $Z_f(n)$ to the set of quadratics $\cQ_f$ as follows. Let $\cQ_f[n] = \{ [u,r,v]\in\cQ_f: v=n\}$; plainly the group $\Gamma_\infty=\{\begin{pmatrix} a & b \\ c & d \end{pmatrix} \in\Gamma: c=0\}$ leaves invariant the set $\cQ_f[n]$. On writing $[u,r,v] = [\frac{\eta^2-\Delta}{4n\alpha},\eta,n\alpha]$, one readily shows that the set $Z'_f(n)$, hence also the set $Z_f(n)$, is in bijective correspondence with the set $\Gamma_\infty \backslash \cQ_f[n\alpha]$. Via this bijection, we derive the following identity.

\begin{lemma}\label{lem:pre-Weyl-Poincare} We have
\begin{equation}\label{eq:pre-Weyl-Poincare}
\cW_{h,g}(x,N) = \sum_{j=1}^{h_\Delta} 
	\sum_{\substack{\xi\in \Gamma_\infty \backslash \Gamma / \Gamma^{(j)} \\ 
		\uvj(\xi) \equiv 0 \, (\rm{mod}\, N\alpha)}} 
	g\(\frac{\uvj(\xi)}{\alpha}\) \,\, e\(\frac{h (\urj(\xi)-\beta)}{2 \uvj(\xi)}\) .
\end{equation}
\end{lemma}

\subsection{Congruence groups and divisibility} 


We seek a group-theoretic interpretation of the arithmetic restriction appearing in the expression \eqref{eq:pre-Weyl-Poincare} of $\cW_{h,g}(x,N)$. From this moment on, we let $\Gamma'=\Gamma_0(N\alpha) \subset \Gamma = \Gamma_0(\alpha)$, and put $N_1=N\alpha$ so that $\Gamma'=\Gamma_0(N_1)$.
Our interpretation starts with the observation that for an arbitrary element $\xi \in \Gamma$, the divisibility condition $N_1 | \uvj(\xi)$ depends only on the coset $\Gamma'\xi \in \Gamma'\backslash\Gamma$.

\begin{lemma}\label{lem:divisibility}
Let 
$\xi = \begin{pmatrix} p & q \\ r & s \end{pmatrix} \in \Gamma'$, 
$\eta = \begin{pmatrix} a & b \\ c & d \end{pmatrix}\in\Gamma$, 
and $1 \leq j \leq h_\Delta$. 
Then $\uvj(\xi\eta) \equiv s^2 \uvj(\eta)  \, (\rm{mod}\, N_1).$
In particular, $N_1 | \uvj(\eta)$ if and only if $N_1 | \uvj(\xi\eta)$.
\end{lemma}

\begin{proof}
By definition we have $N_1|r$ and $\gcd(N_1,s)=1$. It follows from \eqref{eq:group-action-v} that
$$
\uvj(\xi\eta) = u_j \uc(\xi\eta)^2 + r_j \uc(\xi\eta) \ud(\xi\eta) + v_j \ud(\xi\eta)^2.
$$
It is evident that $\uc(\xi\eta) \equiv sc \, (\rm{mod}\, N_1)$ and $\ud(\xi\eta) \equiv sd \, (\rm{mod}\, N_1)$, and hence 
$$
\uvj(\xi\eta) \equiv s^2(u_jc^2 + r_jcd + v_jd^2) = s^2 \uvj(\eta)  \, (\rm{mod}\, N_1).
$$
The lemma follows.
\end{proof}

Consider the set of cosets
\begin{equation}\label{eq:divisibility-set}
V_j(N_1;\Gamma'\backslash\Gamma) = \{ \Gamma'\eta \in \Gamma'\backslash\Gamma : \uvj(\eta) \equiv 0 \, (\rm{mod}\, N_1) \};
\end{equation}
this set is well-defined by Lemma \ref{lem:divisibility}. An element in $V_j(N_1;\Gamma'\backslash\Gamma)$ is called a \emph{divisible coset}.
Each divisible coset has a natural action of the corresponding isotropy group on the right as follows. 
If $\Gamma'\eta\in V_j(N_1;\Gamma'\backslash\Gamma)$, $\gamma'\in\Gamma'$, and $\gamma_j\in \Gamma^{(j)}$, set 
\begin{equation}\label{eq:divisibility-coset-group-action}
( \gamma'\eta ) \cdot \gamma_j = \gamma' (\eta\gamma_j\eta^{-1}) \eta.
\end{equation}
That this action is well-defined is the content of the following lemma.

\begin{lemma}
We have $\eta \gamma_j \eta^{-1} \in \Gamma'$.
\end{lemma}

\begin{proof}
Write $\eta = \begin{pmatrix} a & b \\ c & d \end{pmatrix} \in \Gamma$, so $\eta^{-1} = \begin{pmatrix} d & -b \\ -c & a \end{pmatrix}$.
By \eqref{eq:isotropy-group1}, write
$\gamma_j = \begin{pmatrix} \frac{\tau+\upsilon r_j}{2} & -\upsilon u_j \\ \upsilon v_j & \frac{\tau-\upsilon r_j}{2} \end{pmatrix} \in \Gamma^{(j)}$ where $\tau,\upsilon\in\bZ, \tau^2 - \Delta \upsilon^2 = 4$. 
A straightforward matrix multiplication shows that
$$
\eta \gamma_j \eta^{-1} = \begin{pmatrix} * & * \\ \upsilon \uvj(\eta) & * \end{pmatrix} 
$$
By assumption $\Gamma'\eta\in V_j(N_1;\Gamma'\backslash\Gamma)$, we have $N_1 | \uvj(\eta)$. The lemma follows.
\end{proof}

\begin{definition}\label{defn:pre-unitypartition-Poincare-series}
For $\Gamma'\eta \in V_j(N_1;\Gamma'\backslash\Gamma)$, let 
\begin{equation}\label{eq:pre-unitypartition-Poincare-series}
P_j(\Gamma'\eta;x,N) 
	= \sum_{\xi\in \Gamma_\infty \backslash \Gamma' \eta / \Gamma^{(j)}} 
	e\(\frac{h (\urj(\xi)-\beta)}{2 \uvj(\xi)}\) \,\, g\(\frac{\uvj(\xi)}{\alpha}\) .
\end{equation}
\end{definition}

\begin{corollary}\label{cor:pre-Weyl-Poincare} 
We have
\begin{equation}\label{eq:pre-unitypartition-Weyl-Poincare}
\cW_{h,g}(x,N) = \sum_{j=1}^{h_\Delta}  \,\,\,
	\sum_{\Gamma'\eta \in V_j(N_1;\Gamma'\backslash\Gamma)} P_j(\Gamma'\eta;x,N).
\end{equation}
\end{corollary}

\begin{proof}
This is an immediate consequence of Lemma \ref{lem:pre-Weyl-Poincare} and \eqref{eq:pre-unitypartition-Poincare-series}.
\end{proof}


\begin{lemma}\label{lem:number-of-divisible-cosets}
For $1\leq j \leq h_\Delta$, we have
$$
|V_j(N_1; \Gamma'\backslash\Gamma)| \ll \tau(N).
$$
\end{lemma}

\begin{proof}
First note that $V_j(N_1; \Gamma'\backslash\Gamma) \subset \Gamma'\backslash\Gamma$.
Define the map 
$$\varphi: \Gamma'\backslash\Gamma \to \Gamma'\backslash\bP^1(\bQ),  
\varphi(\Gamma'\eta)= [\eta\cdot\infty]$$ 
where for $\fa\in\bQ$ we let $[\fa]$ denote the orbit of $\fa$ in $\Gamma'\backslash\bP^1(\bQ)$. 
We want to bound $|V_j(N_1; \Gamma'\backslash\Gamma)|$ by considering the image $\varphi(V_j(N_1; \Gamma'\backslash\Gamma))$ and fibers of $\varphi$.

By Lemma \ref{lem:cusps-representatives-width}, a complete set of inequivalent cusps of $\Gamma'$ is given by 
$$ \Big\{ \frac{\mu}{\nu}: \mu,\nu \in \bZ_{>0}, \nu | N_1, \gcd(\mu,\nu)=1, \mu \textrm{ mod } \, \gcd(\nu, N_1/\nu) \Big\}.$$
It is apparent that for each $\eta\in\Gamma$, the point $\eta\cdot\infty\in\bP^1(\bQ)$ is a cusp for $\Gamma'$.  We infer that each coset in $\Gamma'\backslash\Gamma$ can be written as $\Gamma'\eta$ where $\eta = \begin{pmatrix} a & b \\ c & d \end{pmatrix} \in \Gamma$ satisfies $c | N_1$ and $1 \leq a < c$. Now the defining condition \eqref{eq:divisibility-set} 
for $\Gamma'\eta \in V_j(N_1; \Gamma'\backslash\Gamma)$ is
$N_1 | (u_jc^2 + r_jcd + v_jd^2).$
It follows that $c | (u_jc^2 + r_jcd + v_jd^2)$, and so $c | v_j$. In particular $c=O(1)$ and hence $a=O(1)$. We infer that the set $\varphi(V_j(N_1; \Gamma'\backslash\Gamma))$, being a subset of 
$\{[a/c] \in \Gamma'\backslash\bP^1(\bQ): a=O(1),c=O(1)\}$, has size $O(1)$.

Now fix a divisible coset $\Gamma'\eta_0 \in V_j(N_1; \Gamma'\backslash\Gamma)$, where 
$\eta_0 = \begin{pmatrix} a_0 & b_0 \\ c_0 & d_0 \end{pmatrix} \in \Gamma$ satisfies $c_0 | N_1$ and $1 \leq a_0 < c_0$. By the discussion in the previous paragraph, we may assume $c_0=O(1)$.
We claim that 
$\varphi^{-1}(\varphi(\Gamma'\eta_0)) \cap V_j(N_1; \Gamma'\backslash\Gamma)$ 
has size $O(\tau(N))$, whence the lemma follows. The proof of the lemma is reduced to verifying the claim.

Let $\Gamma'\eta \in V_j(N_1; \Gamma'\backslash\Gamma)$ 
be such that $\varphi(\Gamma'\eta) = \varphi(\Gamma'\eta_0)$. This means there exists $\gamma\in\Gamma'$ such that $\eta_0\cdot\infty = \gamma\eta \cdot \infty$; in other words there exists 
$\gamma_\infty(k) = \begin{pmatrix} 1 & k \\ 0 & 1 \end{pmatrix} \in \Gamma_\infty$ satisfying $\gamma\eta=\eta_0\gamma_\infty(k)$. Hence $\Gamma'\eta = \Gamma' \eta_0\gamma_\infty(k)$. We now observe that if $k\equiv k' \, (\rm{mod}\, N)$, then 
$\Gamma'\eta_0\gamma_\infty(k) = \Gamma'\eta_0\gamma_\infty(k')$. To show this observation, first write $k'=k+fN$ for $f\in\bZ$, and then note that
$$
\eta_0\gamma_\infty(k') = \eta_0\gamma_\infty(fN) \gamma_\infty(k) 
	= (\eta_0\gamma_\infty(fN)\eta_0^{-1}) \cdot (\eta_0\gamma_\infty(k)) 
$$
and that $\eta_0\gamma_\infty(fN)\eta_0^{-1} = 
\begin{pmatrix} * & * \\ -c_0fN & * \end{pmatrix} \in \Gamma'$. This observation implies that the size of the set $\varphi^{-1}(\varphi(\Gamma'\eta_0)) \cap V_j(N_1; \Gamma'\backslash\Gamma)$ does not exceed the number of $0\leq k < N$ such that $\Gamma'\eta_0\gamma_\infty(k) \in V_j(N_1; \Gamma'\backslash\Gamma)$. Since both $\Gamma'\eta_0$ and $\Gamma'\eta_0\gamma_\infty(k)$ are divisible cosets, we have 
\begin{align*}
N_1 &| u_jc_0^2 + r_jc_0d_0 + v_jd_0^2 , \\
N_1 &| u_jc_0^2 + r_jc_0(d_0+c_0k) + v_j(d_0+c_0k)^2 .
\end{align*}
It follows that $N_1 | c_0k (r_jc_0 + v_jc_0k + 2v_jd_0)$. We note that $v_j\neq 0$ since otherwise $\Delta=r_j^2-4u_jv_j$ would be a square. Write $N_2=\gcd(N_1,k)$, $N_1=N_2N_2'$ and $k=N_2k'$. The number of $N_2$ is $O(\tau(N))$, and for each such $N_2$ the number of $k$ which is divisible by $N_2$ and which satisfies $N_1 | c_0k (r_jc_0 + v_jc_0k + 2v_jd_0)$ is $O(1)$. Thus the number of $0\leq k < N$ such that $\Gamma'\eta_0\gamma_\infty(k) \in V_j(N_1; \Gamma'\backslash\Gamma)$
is $O(\tau(N))$. The proof is complete.
\end{proof}

\begin{lemma}\label{lem:divisible-width-asymptotic}
Suppose that $1\leq j \leq h_\Delta$.
If $\Gamma'\eta \in V_j(N_1; \Gamma'\backslash\Gamma)$, then 
${\rm width}_{\Gamma'}(\Gamma'\eta) \asymp N$.
\end{lemma}

\begin{proof}
From the proof of Lemma \ref{lem:number-of-divisible-cosets}, we see that a divisible coset in $V_j(N_1; \Gamma'\backslash\Gamma)$ can be written as $\Gamma'\eta$ with $\eta = \begin{pmatrix} a & b \\ c & d \end{pmatrix} \in \Gamma$ satisfying $c | N_1$, $1 \leq a < c$, $a=O(1)$, and $c=O(1)$.
It follows from Lemma \ref{lem:cusps-representatives-width} that
$$
{\rm width}_{\Gamma'}(\Gamma'\eta) = \frac{N_1}{\gcd(N_1,c^2)} \asymp N_1 \asymp N.
$$
The lemma follows.
\end{proof}

\subsection{Partition of unity}\label{subsection:partition-of-unity}

To deal with the groups of automorphs $\Gamma^{(j)}$, we make use of an ingenious device due to T\'oth \cite{Toth}. 

\begin{definition}
Let $T \in \SL_2(\bZ)$ be such that $T$ has exactly two distinct nonzero real fixed points $t_0<t_1$. A function $p:\bR\to\bR_{\geq 0}$ is called a \emph{$T$-function} if it satisfies the following conditions:
\begin{enumerate}
\item[${\rm (i)}$] $p$ is smooth and compactly supported; 
\item[${\rm (ii)}$] $\|p^{(j)}\|_\infty \ll_j 1$ for all $j\geq 0$;
\item[${\rm (iii)}$] there exist $0<r_1<r_2$ such that $p(x)\neq 0$ implies $r_1<|x|<r_2$;
\item[${\rm (iv)}$] For every $x\in \bQ$, one has $\sum_{n\in\bZ} p(T^n(x)) = 1$ where $T^n(x)=T\circ T\circ \cdots T(x)$ denotes the $n$-fold composition of $T$.
\end{enumerate}
\end{definition}

\begin{lemma}\label{lem:Toth-function}
Suppose that $T\in\SL_2(\bZ)$ has exactly two distinct nonzero real fixed points $t_0<t_1$. Then there exists a $T$-function.
\end{lemma}

\begin{proof} 
Let us construct the function $p:(t_0,t_1)\to\bR_{\geq 0}$; the construction of $p$ outside of $(t_0,t_1)$ can be done by a similar argument.
Since $t_0<t_1$ are nonzero reals, we either have $t_1>0$ or $t_0<0$. Without loss of generality we assume that $t_1>0$; the other case can be treated similarly.

We first observe that if $t_0<t<t_1$, then $\lim_{n\to +\infty}T^n(t)=t_1$ and $\lim_{n\to -\infty}T^n(t)=t_0$. With this observation, we can begin our construction.

In the interval $(t_0,t_1)$ we select several numbers as follows.
We choose an arbitrary real number $u_1$ such that $\max(0,t_0) < u_1 < t_1$ and let $u'_1=T(u_1)$, $u''_1=T(u'_1)$. We choose a real number $u_2$ such that $u_1<u_2<u'_1$ and let $u'_2=T(u_2)$.
Now let $p_0:\bR\to\bR_{\geq 0}$ be a smooth and compactly supported function satisfying the following conditions: 
\begin{itemize}
\item $p_0=1$ on $[u_2,u'_2]$; 
\item $p_0$ is supported in the set $(u_1,u''_1)$; 
\item $\|p_0^{(j)}\|_\infty \ll_j 1$ for all $j\geq 0$.
\end{itemize}
By construction we have $\sum_{n\in\bZ} p_0(T^n(t)) > 0$ for every $t_0<t<t_1$. We now define for $t_0<t<t_1$ the function 
$$ p(t) = \frac{p_0(t)}{\sum_{n\in\bZ} p_0(T^n(t))} .$$
One readily verifies that the function $p$ satisfies the conditions of a $T$-function on the interval $(t_0,t_1)$.
\end{proof}

Let us now apply the above construction to the following situation. 

\begin{proposition}\label{prop:Toth-pou}
Let $\Gamma=\Gamma_0(\alpha)$ and let $q=[u,r,v]$ be a binary quadratic form with integer coefficients of positive non-square discriminant $r^2-4uv=\Delta$ such that $\alpha | v$. 
Let $\Gamma^{q} = \{\gamma\in\Gamma: \gamma\cdot q = q\}.$

There exists a function $\psi_q: \Gamma \to \bR_{\geq 0}$ satisfying the following properties: 
\begin{itemize}
\item $\psi_q(\gamma_\infty\xi)=\psi_q(\xi)$ for all $\xi\in\Gamma,\gamma_\infty\in\Gamma_\infty$; 
\item if $\xi\in\Gamma$, then $\sum_{\gamma\in\Gamma^q} \psi_q(\xi\gamma)=1$.
\item there exist $0<r_1<r_2$ such that 
if $\xi=\begin{pmatrix} a & b \\ c & d \end{pmatrix} \in \Gamma$ satisfies $\psi_q(\xi) \neq 0$, then $r_1< \frac{c}{d} < r_2$.
\end{itemize}
\end{proposition}

\begin{proof}
Let $(\tau_0,\upsilon_0)$ be the fundamental solution of the Pell equation $\tau^2-\Delta \upsilon^2=4$.
Set $T_0 = \begin{pmatrix} \frac{\tau_0+\upsilon_0 r}{2} & -\upsilon_0 u \\ 
			\upsilon_0 v & \frac{\tau_0-\upsilon_0 r}{2} \end{pmatrix}$; 
let $T_1$ denote the transpose of $T_0$.
We have
$$\Gamma^{q} = \{\gamma\in\Gamma: \gamma\cdot q = q\}
	= \bigg\{ \pm T_0^n : n \in \bZ \bigg\} . $$

Let $p_1$ be a $T_1$-function. 
We define the function $\psi_q: \Gamma \to \bR_{\geq 0}$ as follows. If 
$\xi=\begin{pmatrix} a & b \\ c & d \end{pmatrix} \in \Gamma$, set $\psi_q(\xi) = p_1\(\frac{c}{d}\)$. 
The required properties of $\psi_q$ follow from the corresponding properties of the $T_1$-function $p_1$.
\end{proof}

With the partition of unity construction, we can unfold the right hand side of \eqref{eq:pre-unitypartition-Poincare-series} to express $P_j(\Gamma'\eta;x,N)$, where $1 \leq j \leq h_\Delta$ and $\Gamma'\eta\in V_j(N_1;\Gamma'\backslash\Gamma)$, as a Poincar\'e series. Let us abbreviate $\psi_j=\psi_{q_j}$ for $1 \leq j \leq h_\Delta$.
\begin{corollary}\label{cor:pre-Poincare-series}
For $1 \leq j \leq h_\Delta$ and $\Gamma'\eta \in V_j(N_1;\Gamma'\backslash\Gamma)$, we have 
\begin{equation}\label{eq:pre-Poincare-series}
P_j(\Gamma'\eta;x,N) 
	= \sum_{\xi\in \Gamma_\infty \backslash \Gamma' \eta } 
	e\(\frac{h (\urj(\xi)-\beta)}{2 \uvj(\xi)}\) \,\, g\(\frac{\uvj(\xi)}{\alpha}\) \psi_j(\xi).
\end{equation}
\end{corollary}

We invoke Hooley's identity to approximate $\cW_{h,g}(x,N)$ by Poincar\'e series of simpler form.

\begin{corollary}\label{cor:Poincare-series}
For $1 \leq j \leq h_\Delta$ and $\Gamma'\eta \in \Gamma'\backslash\Gamma$, let 
\begin{equation}\label{eq:Poincare-series-definition}
Q_j(\Gamma'\eta;x,N) 
	= \sum_{\xi = \begin{pmatrix} a & b \\ c & d \end{pmatrix} 
	\in \Gamma_\infty \backslash \Gamma' \eta } 
	e\(\frac{ah}{c}\) \,\, g\(\frac{\uvj(\xi)}{\alpha}\) \psi_j(\xi).
\end{equation}
If $\Gamma'\eta \in V_j(N_1;\Gamma'\backslash\Gamma)$, then
\begin{equation}\label{eq:Weyl-Poincare-PQ}
P_j(\Gamma'\eta;x,N)  =  Q_j(\Gamma'\eta;x,N)  +  O(h).
\end{equation}
Consequentially,
\begin{equation}\label{eq:Weyl-Poincare}
\cW_{h,g}(x,N) = \sum_{j=1}^{h_\Delta}  \,\,\,
	\sum_{\Gamma'\eta \in V_j(N_1;\Gamma'\backslash\Gamma)} Q_j(\Gamma'\eta;x,N)
	+ O(h\tau(N)).
\end{equation}
\end{corollary}

\begin{proof}
We consider $\xi = \begin{pmatrix} a & b \\ c & d \end{pmatrix} \in \Gamma' \eta$ for which $ g\(\frac{\uvj(\xi)}{\alpha}\) \psi_j(\xi) \neq 0$. By Proposition \ref{prop:Toth-pou}, we have $c \asymp \sqrt{x} \asymp d$ and $\uvj(\xi) \asymp x$.
It then follows from Lemma \ref{lem:Hooley-identity} that
$$
\frac{h(\urj(\xi)-\beta)}{2\uvj(\xi)} 
	= \frac{ah}{c} - \frac{h(r_jc+2v_jd)}{2c \uvj(\xi)} - \frac{h\beta}{2\uvj(\xi)}
	= \frac{ah}{c} + O\(\frac{h}{x}\).
$$
Therefore $e\(\frac{h(\urj(\xi)-\beta)}{2\uvj(\xi)}\) = e\( \frac{ah}{c} \) + O\(\frac{h}{x}\).$
We deduce that
\begin{align*}
P_j(\Gamma'\eta;x,N) 
	&= \sum_{\xi = \begin{pmatrix} a & b \\ c & d \end{pmatrix} 
	\in \Gamma_\infty \backslash \Gamma' \eta } 
	e\(\frac{ah}{c}\) \,\, g\(\frac{\uvj(\xi)}{\alpha}\) \psi_j(\xi)
	+ O\( \sum_{c\asymp \sqrt{x}, d\asymp \sqrt{x}} \frac{h}{x} \) \\
	&= Q_j(\Gamma'\eta;x,N) + O(h).
\end{align*}
This proves \eqref{eq:Weyl-Poincare-PQ}. On applying Lemma \ref{lem:number-of-divisible-cosets}, we deduce the approximation \eqref{eq:Weyl-Poincare}. 
\end{proof}

\subsection{Spectral analysis of Poincar\'e series}\label{subsection:Poincare-Kloosterman}


Recall that $\Gamma'=\Gamma_0(N_1)=\Gamma_0(N\alpha)\subset \Gamma=\Gamma_0(\alpha)$.
Suppose that $\eta \in \SL_2(\bZ)$ and that $\fa=\eta\cdot\infty=\frac{\mu}{\nu}\in\bQ$ is a cusp for $\Gamma_0(N_1)$ with $\gcd(\mu,\nu)=1$ and $\nu | N_1$. Put $N_1=\nu N'$, $\iota=\gcd(\nu,N')$, $\nu=\iota\nu'$, and $N'=\iota N''$. Let $\sigma_\fa = \eta \tau_\fa$, where $\tau_\fa=\begin{pmatrix} \sqrt{N''} & 0 \\ 0 & 1 / \sqrt{N''} \end{pmatrix}$, be the scaling matrix chosen for the cusp $\fa$.
Let
\begin{align}
g_j(c,\kappa,y) 
	&= g\(\frac{q_j(c,y)}{\alpha}\)  \psi_j(c,y)  e\(\frac{-y\kappa}{cN''}\)  \label{eq:Fourier-integrand} \\
G_j(c,\kappa) 
	&= \int_{-\infty}^{\infty} g\(\frac{q_j(c,y)}{\alpha}\)  \psi_j(c,y)  e\(\frac{-\kappa y}{cN''}\)   \,{\rm dy}
	= \int_{-\infty}^{\infty} g_j(c,\kappa,y)  \,{\rm dy} . \label{eq:Fourier-transform}
\end{align}

\begin{lemma}\label{lem:Poincare-Kloosterman}
We have
\begin{equation}\label{eq:Poincare-Kloosterman}
Q_j(\Gamma'\eta;x,N) =
	\sum_{c \in \cQ(\nu,N')} \sum_{\kappa\in\bZ}
	\frac{S_{\sigma_\infty \sigma_\fa}(h,\kappa;c\sqrt{N''})}{cN''} G_j(c,\kappa).
\end{equation}
\end{lemma}

\begin{proof}
The Kloosterman sums that appear in our analysis of $Q_j(\Gamma'\eta;x,N)$ are related to the cusps $\infty$ and $\fa$ (Definition \ref{defn:Kloosterman-sum}). By Lemma \ref{lem:cusps-representatives-width}, we have ${\rm width}_{\Gamma'}(\Gamma'\eta)=N''$; in other words,
$$
\eta^{-1} \Gamma'_\fa \eta = \pm  \begin{pmatrix} 1 & N'' \bZ \\ 0 & 1 \end{pmatrix}.
$$
It follows that, by \eqref{eq:Poincare-series-definition},
\begin{align*}
Q_j(\Gamma'\eta;x,N) &= 
	\sum_{\xi = \begin{pmatrix} a & b \\ c & d \end{pmatrix} \in \Gamma_\infty \backslash \Gamma' \eta} 
	 e\(\frac{ah}{c}\) \,\, g\(\frac{\uvj(\xi)}{\alpha}\)  \psi_j(\xi) \\
	&= \sum_{\xi = \begin{pmatrix} a & b \\ c & d \end{pmatrix} \in 
			\Gamma_\infty \backslash \Gamma' \eta / \eta^{-1} \Gamma'_\fa \eta } 
		e\(\frac{ah}{c}\) \, \sum_{k \in \bZ} 
	g\(\frac{1}{\alpha}  \uvj\(\xi\cdot \begin{pmatrix} 1 & N''k \\ 0 & 1 \end{pmatrix}\)\)\, \psi_j\(\xi\cdot \begin{pmatrix} 1 & N''k \\ 0 & 1 \end{pmatrix}\) \\
	&= \sum_{\begin{pmatrix} a & b \\ c & d \end{pmatrix} \in 
			\Gamma_\infty \backslash \Gamma' \eta / \eta^{-1} \Gamma'_\fa \eta } 
		e\(\frac{ah}{c}\) \, \sum_{k \in \bZ} 
	g\(\frac{q_j(c,d+cN''k)}{\alpha}\)\, \psi_j(c,d+cN''k) .
\end{align*}

On applying Poisson summation to the $k$-sum, denoting by $\kappa$ the dual variable, and changing variable $y=d+cN''t$, we infer that
\begin{align*}
Q_j(\Gamma'\eta;x,N) 
	&= \sum_{\begin{pmatrix} a & b \\ c & d \end{pmatrix} \in 
			\Gamma_\infty \backslash \Gamma' \eta / \eta^{-1} \Gamma'_\fa \eta } 
		e\(\frac{ah}{c}\) \, \sum_{\kappa \in \bZ} 
	\int_{-\infty}^{\infty} g\(\frac{q_j(c,d+cN''t)}{\alpha}\))\, \psi_j(c,d+cN''t) e(-t\kappa) \, {\rm dt} \\
	&= \sum_{\kappa \in \bZ} \sum_{\begin{pmatrix} a & b \\ c & d \end{pmatrix} \in 
			\Gamma_\infty \backslash \Gamma' \eta / \eta^{-1} \Gamma'_\fa \eta } 
		\frac{1}{cN''} e\(\frac{ah}{c} + \frac{d\kappa}{cN''}\) \,
	\int_{-\infty}^{\infty} g\(\frac{q_j(c,y)}{\alpha}\)\, \psi_j(c,y) e\(\frac{-y\kappa}{cN''}\) \, {\rm dy} \\
	&= \sum_{\kappa \in \bZ} \sum_{\begin{pmatrix} a & b \\ c & d \end{pmatrix} \in 
			\Gamma_\infty \backslash \Gamma' \eta / \eta^{-1} \Gamma'_\fa \eta } 
		\frac{1}{cN''} e\(\frac{ah}{c} + \frac{d\kappa}{cN''}\) \, G_j(c,\kappa) ,
\end{align*}
by \eqref{eq:Fourier-transform}.
It now follows from Lemma \ref{lem:Kloosterman-moduli} that
\begin{align*}
Q_j(\Gamma'\eta;x,N) 
	&= \sum_{\kappa \in \bZ} \sum_{\substack{c \in \cC_{\infty \fa}(\Gamma') / \sqrt{N''} \\ 
		\begin{pmatrix} a & b \\ c & d \end{pmatrix} \in 
			\Gamma_\infty \backslash \Gamma' \eta / \eta^{-1} \Gamma'_\fa \eta }} 
		\frac{1}{cN''} e\(\frac{ah}{c} + \frac{d\kappa}{cN''}\) \, G_j(c,\kappa) \\
	&= \sum_{\kappa\in\bZ} \sum_{c \in \cC_{\infty \fa}(\Gamma') / \sqrt{N''}}
		\frac{S_{\sigma_\infty \sigma_\fa}(h,\kappa;c\sqrt{N''})}{cN''} G_j(c,\kappa) \\
	&= \sum_{\kappa\in\bZ} \sum_{c \in \cQ(\nu,N')}
		\frac{S_{\sigma_\infty \sigma_\fa}(h,\kappa;c\sqrt{N''})}{cN''} G_j(c,\kappa).
\end{align*}
The lemma is proved.
\end{proof}


The following two lemmas provide necessary bounds for our analysis of Poincar\'e series. Their proofs are deferred to Section \ref{subsection:auxiliary-estimate}, so as not to interrupt the flow of our analysis.

\begin{lemma}\label{lem:Poincare-Kloosterman-estimate}
For any integer $R \geq 1$ we have
\begin{equation}\label{eq:Fourier-transform-bound}
G_j(c,\kappa) \ll_R \sqrt{x} \cdot  \(\frac{NY_1}{\kappa}\)^R .
\end{equation}
\end{lemma}

\begin{lemma}\label{lem:Fourier-integrand-derivatives}
For $I,J\in\bN$ we have
\begin{equation}\label{eq:Fourier-integrand-derivatives}
\( \frac{d}{dc} \)^I \( \frac{d}{d\kappa} \)^J  g_j(c,\kappa,y) 
	\ll \( \frac{Y_1}{\sqrt{x}} \)^I \( \frac{1}{N} \)^J (NY_1)^\epsilon.
\end{equation}
\end{lemma}

We are in a position to estimate the Poincar\'e series $Q_j(\Gamma'\eta;x,N)$.

\begin{proposition}\label{prop:Poincare-Kloosterman-restricted}
We have
\begin{equation}\label{eq:Pitt-application-2}
Q_j(\Gamma'\eta;x,N) 
	\ll \( x^\frac34 N^{-\frac12} Y_1^\frac94 \gcd(h,N)^\frac14 
	+ x^\frac12 Y_1^3 	\)  (xNY_1)^\epsilon .
\end{equation}
\end{proposition}

\begin{proof}
By Lemma \ref{lem:divisible-width-asymptotic} we have $\mu,\nu,\iota \asymp 1$ and $N',N'' \asymp N$.
Consider \eqref{eq:Poincare-Kloosterman} of Lemma \ref{lem:Poincare-Kloosterman}.
On applying Lemma \ref{lem:Poincare-Kloosterman-estimate}, we may restrict the $\kappa$-sum in \eqref{eq:Poincare-Kloosterman} to $|\kappa| < (NY_1)^{1+\epsilon}$ modulo an admissible error term. More precisely,
\begin{equation}\label{eq:Poincare-Kloosterman-kappa-restricted}
Q_j(\Gamma'\eta;x,N) = \frac{1}{N''}
	\int_{-\infty}^{\infty} 
	\sum_{c \in \cQ(\nu,N')} \sum_{|\kappa| < (NY_1)^{1+\epsilon}}
	\frac{S_{\sigma_\infty \sigma_\fa}(h,\kappa;c\sqrt{N''})}{c} 
	g_j(c,\kappa,y) \,{\rm dy} + O(N^{-100}) .
\end{equation}

We now apply Theorem \ref{thm:Pitt-estimate} to estimate the integrand of the right hand side of \eqref{eq:Poincare-Kloosterman-kappa-restricted}, verifying the required hypotheses by Lemma \ref{lem:Fourier-integrand-derivatives}, and deduce that

\begin{align}
&\sum_{c \in \cQ(\nu,N')} \sum_{|\kappa| < (NY_1)^{1+\epsilon}}
	\frac{S_{\sigma_\infty \sigma_\fa}(h,\kappa;c\sqrt{N''})}{c} g_j(c,\kappa,y)
	\nonumber \\
	&\ll (NY_1)^{\frac12} \bigg\{     
		 \( (NY_1)^\frac14 N^{-\frac14} + 1 \) x^\frac14 Y_1^\frac32 \gcd(h,N)^\frac14  
		+ N^{\frac12}\( Y_1^\frac34 + (NY_1)^\frac12 N^{-\frac12} \) Y_1^\frac74
		\bigg\} 
	(xNY_1)^\epsilon  
	\nonumber \\
	&\ll \( x^\frac14 N^\frac12 Y_1^\frac94 \gcd(h,N)^\frac14 
		+ N Y_1^3	\)  (xNY_1)^\epsilon .   
	\label{eq:Pitt-application-1}
\end{align}
Thus
\begin{align*}
Q_j(\Gamma'\eta;x,N) 
	&\ll \frac{\sqrt{x}}{N} 
		\( x^\frac14 N^\frac12 Y_1^\frac94 \gcd(h,N)^\frac14 
		 + N Y_1^3	\)  (xNY_1)^\epsilon  	\\
	&\ll \( x^\frac34 N^{-\frac12} Y_1^\frac94 \gcd(h,N)^\frac14 
	+ x^\frac12 Y_1^3 	\)  (xNY_1)^\epsilon .
\end{align*}
The proposition is proved. 
\end{proof}

\begin{corollary}\label{cor:smooth-linear-form}
We have
\begin{equation}\label{eq:smooth-linear-form-estimate}
\cW_{h,g}(x,N) 
	\ll  \( x^\frac34 N^{-\frac12} Y_1^\frac94 \gcd(h,N)^\frac14
	+ x^\frac12 Y_1^3  \)  
		(xNY_1)^\epsilon   +   h\tau(N) .
\end{equation}
\end{corollary}

\begin{proof}
In \eqref{eq:Weyl-Poincare}, we apply Lemma \ref{lem:number-of-divisible-cosets} and Proposition \ref{prop:Poincare-Kloosterman-restricted} to obtain the bound \eqref{eq:smooth-linear-form-estimate}.
\end{proof}

We are in a position to prove our main theorem.

\begin{proof}[Proof of Theorem \ref{thm:discrete-linear-form-estimate}]
It follows from Lemma \ref{lem:approximate-linear-form} and Corollary \ref{cor:smooth-linear-form}
that
$$
\cW_{h}(x,N) \ll  \frac{x^{1+\epsilon}}{NY_1}  
			+ \(  
				x^\frac34 N^{-\frac12} Y_1^\frac94 \gcd(h,N)^\frac14
				+ x^\frac12 Y_1^3  \)  
			(xNY_1)^\epsilon   
			+   h\tau(N) .
$$
On choosing $Y_1 = \(\frac{x}{N^2h}\)^{\frac{1}{13}}$, we deduce that
$$
\cW_{h}(x,N) \ll  x^{\frac{12}{13}+\epsilon} N^{-\frac{11}{13}} h^{\frac{1}{13}}  +  h\tau(N) .
$$
The theorem is proved.
\end{proof}

\subsection{Proofs of auxiliary estimates}\label{subsection:auxiliary-estimate}

Recall that $g$ is a smooth weight function supported on the interval $[x-\frac{x}{Y_1},2x+\frac{x}{Y_1}]$ and that its derivatives are bounded by $\| g^{(j)} \|_\infty \ll_j \(\frac{Y_1}{x}\)^j$ for all $j\in\bN$.

\begin{lemma}\label{lem:derivative-bound-1}
Suppose that $q(X,Y)=\alpha X^2 + \beta XY + \gamma Y^2$ and that $c,y \ll \sqrt{x}$. 
If $R\geq 1$, then 
$$
\(\frac{d}{dy}\)^R g(q(c,y)) \ll_{\alpha,\beta\gamma} \(\frac{Y_1}{\sqrt{x}}\)^R.
$$
\end{lemma}

\begin{proof}
Put $\partial=\frac{d}{dy}$. 
We have $\partial q(c,y) = \beta c + 2\gamma y$ and $\partial^2 q(c,y) = 2\gamma$. If $c,y \ll \sqrt{x}$, then $\partial q(c,y) \ll \sqrt{x}$ and $\partial^2 q(c,y) \ll 1$.
By induction we see that there are integer constants $r_j$ so that
$$
\partial^R g(q(c,y)) = \sum_{0\leq j \leq \frac{R}{2}} 
			r_j f^{(R-j)}(q(c,y)) \cdot \(\partial q(c,y)\)^{R-2j} \cdot \(\partial^2 q(c,y)\)^{j}.
$$
Therefore, for  $c,y \ll \sqrt{x}$ we have
$$
\partial^R g(q(c,y)) \ll \max_{0\leq j \leq \frac{R}{2}} \(\frac{Y_1}{x}\)^{(R-j)} \cdot (\sqrt{x})^{R-2j} 
	\ll  \(\frac{Y_1}{\sqrt{x}}\)^{R}.
$$
The lemma is proved.
\end{proof}

\begin{lemma}\label{lem:derivative-bound-2}
Let $1\leq j\leq h_\Delta$ and $R\geq 1$. If $c,y \ll \sqrt{x}$ and $|c/y|\asymp 1$, then   
$$
\(\frac{d}{dy}\)^R \psi_j(c,y) \ll \(\frac{1}{\sqrt{x}}\)^R.
$$
\end{lemma}

\begin{proof}
Put $p_{j,1}(x)= p_j(1/x)$, so that
$\psi_j(c,y)=p_j(c/y)=p_{j,1}(y/c)$. Applying Proposition \ref{prop:Toth-pou}, we deduce that
$$
\(\frac{d}{dy}\)^R \psi_j(c,y) = \(\frac{d}{dy}\)^R p_{j,1}(y/c) 
\asymp c^{-R} p_{j,1}^{(R)}(y/c) \asymp \(\frac{1}{\sqrt{x}}\)^R.
$$
The lemma is proved.
\end{proof}

We are in a position to prove Lemma \ref{lem:Poincare-Kloosterman-estimate}.

\begin{proof}[Proof of Lemma \ref{lem:Poincare-Kloosterman-estimate}]
In \eqref{eq:Fourier-transform} we integrate by parts, collecting relevant estimates from Proposition \ref{prop:Toth-pou}, Lemmas \ref{lem:derivative-bound-1} and \ref{lem:derivative-bound-2}, to deduce that, for any integer $R\geq 1$, 
\begin{align*}
G_j(c,\kappa) 
	&= \int_{-\infty}^{\infty} g\(\frac{q_j(c,y)}{\alpha}\)  \psi_j(c,y)  e\(\frac{-\kappa y}{cN''}\)   \,{\rm dy} \\
	&= \int_{-\infty}^{\infty} \(\frac{d}{dy}\)^R \(g\(\frac{q_j(c,y)}{\alpha}\)  \psi_j(c,y)\)  
		\,\, e\(\frac{-\kappa y}{cN''}\) \cdot \(\frac{-2\pi i\kappa}{cN''}\)^{-R}  \,{\rm dy} \\
	&\ll \sqrt{x} \cdot  \(\frac{Y_1}{\sqrt{x}}\)^R \cdot \(\frac{\kappa}{\sqrt{x}N''}\)^{-R}.
\end{align*}
By Lemma \ref{lem:divisible-width-asymptotic} we have $N''\asymp N$. Therefore
$$
G_j(c,\kappa)  \ll \sqrt{x} \cdot  \(\frac{NY_1}{\kappa}\)^R .
$$
The lemma is proved.
\end{proof}

\begin{lemma}\label{lem:derivative-bound-3}
Suppose that $c,y \asymp \sqrt{x}$ and that $\kappa < (NY_1)^{1+\epsilon}$. 
If $R\geq 1$, then 
$$
\(\frac{d}{dc}\)^R e\(\frac{-\kappa y}{cN''}\) \ll \(\frac{Y_1}{\sqrt{x}}\)^R (NY_1)^{\epsilon R}.
$$
\end{lemma}

\begin{proof}
Put $K=-\frac{\kappa y}{N''}$. Then $K \ll \sqrt{x}Y_1 (NY_1)^\epsilon$, by assumptions and by Proposition \ref{prop:Toth-pou}.
By induction we see that there are integer constants $r_j$ so that
$$
\(\frac{d}{dc}\)^R e\(\frac{K}{c}\) = \frac{e\(\frac{K}{c}\)}{c^R} \cdot 
	\bigg\{  \sum_{0\leq j \leq R} r_j \(\frac{K}{c}\)^j  \bigg\}.
$$
It follows that
$$
\(\frac{d}{dc}\)^R e\(\frac{K}{c}\)  \ll \(\frac{Y_1}{\sqrt{x}}\)^R \max_{0\leq j \leq R} (NY_1)^{\epsilon j}.
$$
The lemma follows.
\end{proof}

\begin{proof}[Proof of Lemma \ref{lem:Fourier-integrand-derivatives}]
Recall from \eqref{eq:Fourier-integrand} that
$ g_j(c,\kappa,y) = g\(\frac{q_j(c,y)}{\alpha}\)  \psi_j(c,y)  e\(\frac{-y\kappa}{cN''}\)$. It suffices to show that, for $I,J\in\bN$, 
\begin{align}
\( \frac{d}{dc} \)^I  g_j(c,\kappa,y) &\ll \( \frac{Y_1}{\sqrt{x}} \)^I (NY_1)^\epsilon ,
\label{eq:Fourier-integrand-derivatives-1} \\
\( \frac{d}{d\kappa} \)^J  g_j(c,\kappa,y) &\ll \( \frac{1}{N} \)^J .
\label{eq:Fourier-integrand-derivatives-2}
\end{align}
The bound \eqref{eq:Fourier-integrand-derivatives-2} is immediate on noting that
$$
\( \frac{d}{d\kappa} \)^J  g_j(c,\kappa,y) =  g_j(c,\kappa,y) \(\frac{-2\pi i y}{cN''}\)^J
	\ll \(\frac{\sqrt{x}}{\sqrt{x}N}\)^J = \( \frac{1}{N} \)^J .
$$
We now show \eqref{eq:Fourier-integrand-derivatives-1}. Let $R$ be an arbitrary positive integer. 
By symmetry, we infer from Lemma \ref{lem:derivative-bound-1} that
$\(\frac{d}{dc}\)^R g\(\frac{q_j(c,y)}{\alpha}\) \ll \(\frac{Y_1}{\sqrt{x}}\)^R$, and from Lemma \ref{lem:derivative-bound-2} that
$\(\frac{d}{dc}\)^R \psi_j(c,y) \ll \(\frac{1}{\sqrt{x}}\)^R.$
By Lemma \ref{lem:derivative-bound-3}, we have
$$
\(\frac{d}{dc}\)^R e\(\frac{-\kappa y}{cN''}\) \ll \(\frac{Y_1}{\sqrt{x}}\)^R (NY_1)^{\epsilon}.
$$
The bound \eqref{eq:Fourier-integrand-derivatives-1} thus follows. 
\end{proof}

\subsection*{Acknowledgements}
This work started during my postdoctoral fellowship at the Vietnam Institute for Advanced Study in Mathematics (VIASM) and was completed during my next visit there. I am grateful to VIASM for financial support and hospitality. I would like to thank Hung M.~Bui and Rizwanur Khan for helpful conversations. 



\end{document}